\documentclass{amsart}
\usepackage{geometry}
\usepackage{amsmath}
\usepackage{amssymb}
\usepackage{ulem}
\usepackage{indentfirst}
\usepackage{enumerate}
\usepackage{epstopdf}
\newtheorem{tw}{Theorem}[section]

\newtheorem{cor}{Corollary}[section]

\newcommand{\m}{\overset{m}}
 \usepackage{graphicx}

  \makeatletter
\@addtoreset{equation}{section}
\makeatother
\makeatletter
\@addtoreset{figure}{section}
\makeatother

\begin{document}

\title{A generalization of the total mean curvature}
\author[K.\ Charytanowicz, W.\ Cie\'slak, W.\ Mozgawa]{Katarzyna CHARYTANOWICZ$\;^{*}$,  Waldemar CIE\'SLAK$\;^{**}$,
 Witold  MOZGAWA$^{***}$}
\address{ \begin{itemize}
   \item[$^*$]
  Institute of Mathematics, Maria Curie-Sk{\l}odowska University, pl.~\ M.~\ Curie-Sk{\l}odowskiej 1, 20-031 Lublin, Poland, k.charytanowicz@gmail.com
     \item[$^{**}$]
  Institute of Mathematics and Information Technology, State School of Higher Education in Che{\l}m,
	ul. Pocztowa 54, 22-100 Che{\l}m, Poland, izacieslak@wp.pl
    \item[$^{***}$]   Institute of Mathematics, Maria Curie-Sk{\l}odowska University, pl.~\ M.~\ Curie-Sk{\l}odowskiej 1, 20-031 Lublin, Poland, mozgawa@poczta.umcs.lublin.pl
    \end{itemize}}

  \subjclass{53A05,  	52A15}
  \keywords{support function,   total mean curvature, Hopf formula}

\begin{abstract}
A special formula for the total mean curvature of an ovaloid is derived. This formula allows us to extend the notion of the mean curvature to the class of boundaries  of  strictly convex   sets.  Moreover,  some integral formula for ovaloids is proved.
\end{abstract}
\maketitle

\indent
\section{Introduction}

In this paper we will consider a class of closed surfaces in $\mathbb R^3$. We denote by $A^*$ the region bounded by a surface $A$. According to (22) in \cite{BZ}, p.\ 140, the classical definition of total mean curvature $M(A)$ of the boundary $A$ of a strictly convex body $A^*$ in $\mathbb R^n$ requires that surface $A$  is $C^2$-smooth. Namely
$$M(A)=\int\limits_{A}\frac{k_1+\ldots +k_{n-1}}{n-1}dA,$$
where $k_1,\ldots, k_{n-1}$ are the principal curvatures of $A$ and $dA$ is the area element of $A$. In the case $n=3$  the mean $\frac{k_1+ k_{2}}{2}$ is denoted by $H$ and called the mean curvature of $A$. Thus throughout the paper we will denote the total mean curvature of $A$ by $\int\limits_AH dA$.

We recall that a $C^2$-smooth closed surface in $\mathbb R^3$  with the positive Gaussian curvature $K=k_1k_2$ is called an ovaloid, see \cite{H}, p.\ 119.   As we stated in the abstract the paper  provides a formula for the total mean curvature of an ovaloid in Euclidean three  dimensional space. While our  method involves a straightforward integration technique it gives an interesting extension of the concept of the mean curvature to non-smooth  objects as the boundaries of   strictly convex bodies in $\mathbb R^3$ and its   intriguing  resemblance of the integrand to the potential function is   promising and attracting.

Formula for the total mean curvature introduced in this paper  coincides with the classical formula on the total mean curvature of  ovaloids. However an  analytic form of our formula allows to use it in a wider scope, namely to the class of closed strictly convex surfaces which lies between the class of boundaries of convex bodies  and the class of ovaloids where one uses the ordinary definition.  This is so since in \cite{BZ}, p. 140, it is shown that the total mean curvature $M(A)$ of the boundary $A$ of a regular convex body  $A^*$ in $\mathbb R^n$ satisfies
\begin{equation}\label{dolar}
  M(A)= nV_{n-2}(A^*),
\end{equation}
where $V_{n-2}(A^*)$ is the cross-sectional measure of order $n-2$ (see \cite{BZ}, Sect. 19.3, and \cite{2}, Sect 30). The cross-sectional measures are the geometric functionals at present known in the literature as quermassintegrals or intrinsic volumes (see \cite{Sc}, Ch. 4 and Ch. 5). As a consequence of \eqref{dolar}, in the general case it is natural to assume just $V_{n-2}(A^*)$, possibly up to a constant factor, as the mean curvature of $A$. This point of view is widely confirmed, for example, by the monograph \cite{Sc}.
In particular, for $n=3$, the total mean curvature of $A$ corresponds to $V_{1}(A^*)$, that is up to constant factor, the breadth of $A$ (see, \cite{2}, (9), p. 72). The formula \eqref{100} of our paper just expresses this fact but it seems   that  resulting  formula \eqref{9} is easier to calculate and generalizes the total mean curvature $\int\limits_AH dA$ defined for $C^2$-smooth ovaloids to  the class   of boundaries of  strictly convex sets in $\mathbb  R^3$. We also note at this moment that an interesting plane counterpart of formula \eqref{9} is introduced and investigated in  our paper \cite{CCM}.

At the end of the paper we are able to calculate our formula \eqref{9} in the case of cube which is not a strictly convex body. These calculations show that the crucial ingredient in this framework is just the notion of pedal surface. By this example we do hope that our formula can be useful in investigations of analytic properties of  convex polyhedrons or convex surfaces.

\section{Main result}
Let us fix an ovaloid $E$ and denote by $\overset{m}P$ the support function of $E$ with respect to a fixed  point $m\in E^*$. We associate with $E$ and $m$ a closed surface $\overset{m}E$ defined as follows
$$w\mapsto \overset{m}P(w)w+m$$
for $w\in S^2$, where $S^2$ denotes the unit sphere. $\overset{m}E$ will be called a pedal surface of the surface $E$ with respect to the point $m$, see \cite{Haz}, p. 113.

We denote the outward normal vector field on the surface $E$ by $N$ and by $K$  its  Gaussian curvature. Let $x(u)$ be a local parametrization of the ovaloid  $E$. In the further part of the paper we use the following standard notations for $\alpha, \beta =1,2\colon$ $x_\alpha$ denote the tangent vectors to the coordinate curves, $g_{\alpha\beta} =\langle x_\alpha, x_\beta\rangle$ are  coefficients of the first fundamental form of $E$, $g=\det[g_{\alpha\beta}]$, and if we introduce the following notation $\overset{m}P(N(u))=\overset{m}p(u)$ then
$$\overset{m}p=\langle x-m, N\rangle.$$

Let us fix an ovaloid $E$. We assume that a point $m$ coincides with the origin of the space $\mathbb R^3$ thus we will write $p$ instead $\m p$. Moreover, we parametrize $E$ using the unit sphere $S^2$. We denote by $S_r$ the sphere of radius $r$ with the center at $o$ contained in $E^*$ and by $A_r$ the intersection of the exterior of $S_r$ and $\overset{o} E$. Then we define a mapping
\begin{equation}
\begin{gathered}
F\colon (0,1)\times S^2\to  {A_r}\\
F(t,u)=((1-t)r+tp(u))N(u).
\end{gathered}
\end{equation}
It is easy to see that the Jacobian $F^\prime$ of $F$ is given by the formula
\begin{equation}
F^\prime  = (p-r)((1-t)r+tp)^2K\sqrt g>0.
\end{equation}
Making use of the diffeomorphism $F$ and the Minkowski formula
$$\int\limits_E p K\ dE=\int\limits_E H\  dE,$$
(see \cite{H}, 3), p. 171), we obtain
$$\int\limits_{A_r}\frac{dy}{||y||^2}=\int\limits_E(p-r)K\ dE=\int\limits_E H\  dE-4\pi r,
$$
where $||y||^2=\langle y,y\rangle$ and $dE$ is the surface area form of the ovaloid $E$. On the other hand we have
$$\int\limits_{S_r^*}\frac {dy}{||y||^2}=4\pi r.
$$
Since the integral $\int\limits_EH\ dE$ is   the total mean curvature, so we proved the following theorem.
\begin{tw}
	The total mean curvature $M$ of an ovaloid $E$ is given by the formula
	\begin{equation}\label{4}
	M(E)=\int\limits_{\overset{o}E\raisebox{0.2 em}{\scriptsize *}}\frac {dy}{||y||^2}.
	\end{equation}
\end{tw}

Let the unit sphere $S^2$ be parametrized by
$$w(u)=(\cos u^1\cos u^2, \cos u^1\sin u^2, \sin u^1),$$
for $u=(u^1,u^2)\in U=\left(-\frac \pi 2, \frac \pi 2\right)\times (0, 2\pi).$ Then the volume form is equal to $dS =~-\cos u^1\ du^1\ du^2.$ Note that for any constant vector field $a$ we have
\begin{equation}\label{(20)}
\int\limits_{S^2}\langle a, w\rangle\ dS=0.
\end{equation}
Let $X$ denote the boundary of a fixed strictly
 convex set in $\mathbb R^3$. We assume that the origin $o$ lies in the interior $X^*$ of $X$. We consider two support functions $p_o$ and $p_b$ of X for any $b\in X^*$. It easy to see that
\begin{equation}\label{(21)}
p_b(w)=p_o(w)-\langle ob,w \rangle .
\end{equation}
The conditions \eqref{(20)} and \eqref{(21)} imply that the integral $\int\limits_{S^2}p_b\ dS$ does not depend on $b$. We define the mapping
\begin{equation}
\begin{gathered}
G\colon (0,1)\times U\to \overset{o}X\\
G(t,u)=tp(w(u))w(u),
\end{gathered}
\end{equation}
where $p=p_o$. Since  the  Minkowski support function is differentiable (see \cite{2}, p.29) then we have
\begin{equation}
G^\prime (t,u)=t^2p(w(u))^3\cos u^1.
\end{equation}
Since the  mapping
$$ G\colon (0,1)\times U\to G((0,1)\times U) $$
is a diffeomorphism then we get
\begin{equation}\label{100}
\int\limits_{\overset{o}X\raisebox{0.2 em}{\scriptsize *}}\frac{dy}{||y||^2}=\int\limits_{S^2}p\ dS.
\end{equation}
This formula allows us to extend the notion of the total mean curvature to the class $\mathcal C$ of boundaries of  strictly convex sets in $\mathbb  R^3$. We note that the integral $\int\limits_{\m X\raisebox{0.2 em}{\scriptsize *}}\frac{dy}{||y-m||^2}$ does not depend on $m$. For this reason we can  replace formula \eqref{4} by formula
	\begin{equation}\label{9}
	M(E)=\int\limits_{\m E\raisebox{0.2 em}{\scriptsize *}}\frac{dy}{||y-m||^2},
	\end{equation}
	which is valid for $E$ lying in the class $\mathcal C$ of boundaries of  strictly convex sets in $\mathbb  R^3$.
Thus the function $M\colon  \mathcal C \to \mathbb R$ given by formula \eqref{9} indeed extends the notion of total mean curvature for $C^2$-smooth ovaloids to  $\mathcal C$.

From the above considerations we get an immediate corollary.
\begin{cor}
	If $m$ is a center of the greatest sphere $S_R$ contained in $E^*$, where $E\in \mathcal C$ then the following inequality holds
	\begin{equation}
	M(E)\geq 4\pi R,
	\end{equation}
	and the equality holds if and only if the set $E$ is a sphere.
\end{cor}

Now,  as a by-product we prove some integral formula for ovaloids. Let us define a function
$$p^*\colon \overset{o}E\raisebox{0.2 em}{\scriptsize *}\to \mathbb R$$
as follows -- if $o\not =y\in E^*$ then we denote by $y_E$  the intersection of $E$ and the ray $oy$. We define $p^*(y)$ as the distance between the origin $o$ and  the tangent plane to $E$ at $y_E$ and $p^*(0)=1$. Making use
of the diffeomorphism $F$ we get
$$\int\limits_{A_r}\frac{dy}{(p^*(y)-r)||y||^2}=\int\limits_EK\ dE=4\pi.$$
Hence we get the following theorem.
\begin{tw}
	Let $E$ be an ovaloid. Then the following integral formula holds
	\begin{equation}
	\int\limits_{\overset{o}E\raisebox{0.2 em}{\scriptsize *}}\frac {dy}{p^*(y)||y||^2}=4\pi.
	\end{equation}
\end{tw}

\section{One convex example}
Although there is a substantial body of literature on mean curvature of boundaries of convex bodies we will note that, for example,  the encyclopedia by Santal\'{o}, (\cite{San}, (13.58), p. 226)    gives  a formula for  total mean curvature  for   convex polyhedrons which are a special case of  convex sets. Note also that   Alexander in \cite{ale}  in nontrivial  and deep manner used also this discrete total mean curvature  for closed polyhedral surfaces.

For the convex polyhedrons the edges of which have lengths $a_i$ with corresponding dihedral angles $\alpha_i$  Santalo gives   the following formula for the total mean curvature
\begin{equation}\label{san}
  M=\frac 12 \sum\limits_{i}(\pi-\alpha_i)a_i.
\end{equation}

Note also that  formula (4.23) from  \cite{Sc}
is a particular case of a general connection between intrinsic volumes and generalized curvature
measures which is expressed, for instance, by (4.9) of \cite{Sc}. In $\mathbb R^3$ formula (4.23) implies the above
expression \eqref{san}.

 From the formula \eqref{san} it follows that    the total mean curvature of a cube with sides $2m$ is $6m\pi$, what is in perfect accord with our calculations given below.

  It is interesting that proceeding as for strictly convex bodies in all checked cases we get the correct results. We will illustrate this phenomenon for the case of cube and compare  formula obtained form expression \eqref{9} with a sequence of total mean curvatures of a sequence of ovaloids converging to the cube.
	Consider a cube $C$ of side of length $2m$ and center $(0,0,0)$.	
	\begin{figure}[h]
		\centering
		\includegraphics[height=10cm]{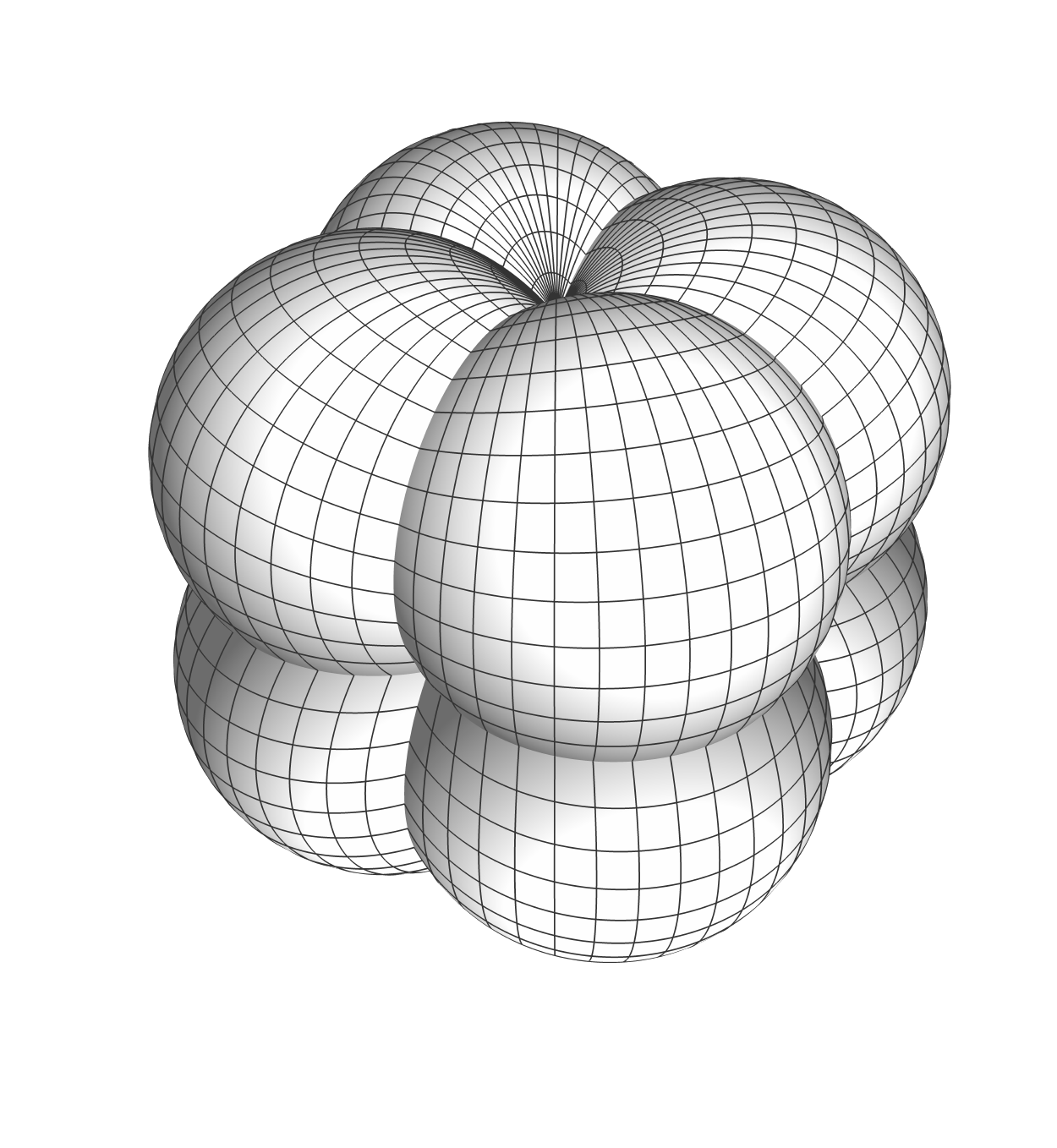}
		\caption{The pedal surface of the cube $C$ with respect to the origin}
	\end{figure}
		
	To calculate the total mean curvature $M$ of the cube $C$ let consider a  part of pedal surface contained in the positive octant denoted by $C_1$. This is the
	part of the sphere $S$ with radius $\frac{\sqrt{3}}{2}m$ centered at the point $\left( \frac{1}{2} m,\frac{1}{2} m,\frac{1}{2} m\right)$. We parametrize it using sphere coordinates centered at the point $\left( \frac{1}{2} m,\frac{1}{2} m,\frac{1}{2} m\right)$ by
	$$w(u)=\left( r(u^1,u^2) \sin u^1 \cos u^2, r(u^1,u^2) \sin u^1 \sin u^2,r(u^1,u^2) \cos u^1 \right), $$
	where $ u=(u^1,u^2)\in U=\left(0, \frac \pi 2\right)\times\left(0, \frac \pi 2\right).$ We determine the term $r(u^1,u^2)$ inserting the coordinates $w(u)$  to the equation of the sphere $S$. Then we get $$r(u^1,u^2)= m \left(  \sin u^1 \sin u^2 +\sin u^1 \cos u^2 + \cos u^1 \right).$$ Hence  we can easily parametrize the interior $\overset{o}{C_1}\raisebox{0.2 em}{\scriptsize *} $ of the part of pedal surface contained in the positive octant by using homothety with center at the origin and a ratio $t \in (0,1)$. Thus 	 we obtain a diffeomorphism
	\begin{equation*}
	\begin{gathered}
	F\colon (0,1) \times U \to \overset{o}{C_1}\raisebox{0.2 em}{\scriptsize *}\\
	F(t,u)=\left( t r(u^1,u^2)  \sin u^1 \cos u^2, t r(u^1,u^2)  \sin u^1 \sin u^2, t r(u^1,u^2)  \cos u^1\right).
	\end{gathered}
	\end{equation*}
	Since the Jacobian  $F'$ of  $F$ is given by
	$$F'=t^2 r^3(u^1,u^2)\sin u^1,$$   by formula \eqref{9}, the total mean curvature $M$ of the cube $C$ is given by
	\begin{equation*}
	M=\iiint\limits_{\overset{o}C\raisebox{0.2 em}{\scriptsize *}} \frac {dx dy dz}{||(x,y,z)||^2}  =
	8 	\int_{0}^{1} \int_{0}^{\frac \pi 2} \int_{0}^{\frac \pi 2}\frac {\left|F' \right|}{ t^2 r^2(u^1,u^2)} du^1 du^2 dt
	\end{equation*}
	and
	\begin{equation*}
	\int_{0}^{1} \int_{0}^{ \frac \pi 2} \int_{0}^{\frac \pi 2}
	\frac {\left|F' \right|}{t^2 r^2(u^1,u^2)} du^1 du^2 dt=\frac{3}{4}m\pi.
	\end{equation*}
	Hence
	\begin{equation*}
	M=8 \cdot \frac{3}{4}m\pi=6 m\pi.
	\end{equation*}
 As we see the above result for the cube coincides with the value obtained from the formula \eqref{san}.

Now we illustrate this result by the numerical calculations for a sequence of strictly convex surfaces  $\displaystyle x^n+y^n+z^n=1 $ with even $n$  tending to the cube with $m=1$  when  $n\to \infty$. The calculations made in {\it Wolfram Mathematica 11.3} with working precision equal to $2000$ confirm the results presented above.

Then the equations  $\displaystyle x^n+y^n+z^n=1 $ with increasing even $n$ describe cube-like surfaces tending to the cube of side of length $2$ and center $(0,0,0)$. The Table \ref{tab:table1} shows the approximate values of the total

\begin{table}[h!]
		
		\caption{Values of $M$ as a function of $n$}
		\begin{center}
			\label{tab:table1}
			\begin{tabular}{c c c c c c }
				\hline
				\textbf{$n$} &$100$ &$200$ &$300$ &$400$ &$500$  \\
				\textbf{$M$} &$18.6792$ & $18.7640$& $18.7928$& $18.8064$& $18.8151$ \\
				\hline
				\\[1pt]
				\hline
				\textbf{$n$} &$600$ &$700$ &$800$ &$900$ &$1000$ \\
				\textbf{$M$} &$18.8208$& $18.8248$ &$18.8280$&  $18.8304$& $18.8312$\\
				\hline
			\end{tabular}
		\end{center}
	\end{table}
\noindent mean curvature $M$ as a function of $n$ for a few values of $n$. It is   seen that following results converge to $6 \pi \approx 18.8496$ as we expected.

\noindent{\bf Acknowledgements:} We would like to thank the referee   for his insightful and substantive comments that helped us to improve the manuscript.

\end{document}